\newtheorem{defi}{Definition}
\newtheorem{theo}{Theorem}
\newtheorem{remark}{Remark}
\newcommand{\proof}{{\em Proof. }}
\newcommand{\qed}{$\hfill\Box$}
\newcommand{\zitat}[4]{\bibitem{#1}{\sc #2}: {\sl #3\/}. #4.\vspace{-0.7em}}
\newcommand{\PG}[2]{\mbox{$\mbox{{\rm PG}}(#1,#2)$}}
\newcommand{\GF}[1]{\mbox{$\mbox{{\rm GF}}(#1)$}}
\newcommand{\abb}[3]{\mbox{$#1\,:\,#2\rightarrow#3$}}

\newcommand{\qu}[1]{\overline{#1}}
\newcommand{\inv}{^{-1}}
\newcommand{\Acal}{{\cal A}}
\newcommand{\Bcal}{{\cal B}}

\newcommand{\Dcal}{{\cal D}}

\newcommand{\Kcal}{{\cal K}}

\newcommand{\Pcal}{{\cal P}}
\newcommand{\Qcal}{{\cal Q}}

\newcommand{\Vcal}{{\cal V}}
\newcommand{\Wcal}{{\cal W}}

\documentstyle[a4,12pt]{article}

\title{A Model of the Witt Design $W_{12}$ based on Quadrics of \PG23%
   \thanks{Research supported by the Austrian FWF, project P12353--MAT.}
   }
\author{Hans Havlicek}
\date{}

\begin{document}
\maketitle

\begin{center}
{\em In memory of H.\ Brauner (1928--1990) on the occasion of his 70th
birthday}
\end{center}

   \begin{abstract}
   An elementary geometric proof for the existence of Witt's $5$--$(12,6,1)$
   design is given.
   \end{abstract}

\noindent Keywords: small Witt design, quadrics.

\section{Introduction}

In the present paper we present a proof for the existence of {\em
Witt's $5$--$(12,6,1)$ design} $W_{12}$. The points of the design
will be all points but one of the projective plane of order three,
the blocks are defined via quadratic equations. Some blocks are
subsets of quadrics, others are sets of points related with quadrics,
e.g., the set of external points of a conic.

Although we shall never make use of it, in the background of our
considerations there will always be the {\em Veronese surface}
$\Vcal$ in \PG53 and a cap $\Kcal$ in \PG53, which is a point model
for Witt's $5$--$(12,6,1)$ design $W_{12}$ \cite{coxe58},
\cite{pell74}. There are various connections between the Veronese
surface $\Vcal$ and the cap $\Kcal$ \cite{havl98}, \cite{havl9y}. By
implementing results from the above--mentioned papers, our proof
could even be shortened. We aim, however, at an elementary proof. In
fact, the prerequisites for reading this article are basic linear
algebra and some properties of quadrics in \PG23.

\section{A planar model of $W_{12}$}

Throughout this paper $F:=\GF3=\{0,1,2\}$ denotes the field with
three elements. The point set $\Pcal(F^3)$ of the projective plane
\PG23 is the set of one--dimensional subspaces of $F^3$. Lines of
\PG23 are considered as sets of points.
   \begin{defi} \label{planar-bloecke}
An incidence structure $(\Wcal,\Bcal,\in)$ with point set $\Wcal$ and
block set $\Bcal$ is given as follows: Fix one point
$U=F(u_0,u_1,u_2)\in\Pcal(F^3)$ and define
   \begin{equation}\label{point-W}
   \Wcal:=\Pcal(F^3)\setminus\{U\}.
   \end{equation}
A subset $b$ of $\Wcal$ is a block, if the subsequent conditions hold
true:
      \begin{enumerate}
      \item $b$ has more than three elements.
      \item There is a non--zero quadratic form \abb{q}{F^3}{F} such that
      $b$ consists of all points $X=F(x_0,x_1,x_2)\in\Wcal$ satisfying
      \begin{equation} \label{planar-koofrei}
      (x_0,x_1,x_2)^q = 2\cdot(u_0,u_1,u_2)^q.
      \end{equation}
      \end{enumerate}
   \end{defi}
Observe that
      \begin{equation}\label{nicht-zweideutig}
      (2\cdot(x_0,x_1,x_2))^q=2^2\cdot(x_0,x_1,x_2)^q=(x_0,x_1,x_2)^q
      \end{equation}
for all $(x_0,x_1,x_2)\in F^3$ and all quadratic forms
\abb{q}{F^3}{F}. Thus (\ref{planar-koofrei}) does not depend on the
choice of vectors representing the points $X$ and $U$, respectively.

It is an easy task to describe all blocks explicitly: Each quadratic
form \abb{q}{F^3}{F} and each $t\in F$ give rise to a point--set of
\PG23 by putting
   \begin{equation}
   \Qcal_t(q):=\{F(x_0,x_1,x_2)\mid (x_0,x_1,x_2)^q=t,\;
   (x_0,x_1,x_2)\neq(0,0,0)\}.
   \end{equation}
By (\ref{nicht-zweideutig}), this definition is unambiguous. Note
that $\Qcal_1(q)=\Qcal_2(2q)$ and $\Qcal_2(q)=\Qcal_1(2q)$.

Up to a change of coordinates and multiplication of $q$ by $2\in F$
there are the following cases for $q\neq 0$ \cite[p.\ 156]{hirs98}:
   \begin{equation}\label{planar-tabelle}
   \mbox{
   \begin{tabular}{|c|c|c|c|}
   \hline
   $(x_0,x_1,x_2)^q$    &$\#\Qcal_0(q)$  &$\#\Qcal_1(q)$  &$\#\Qcal_2(q)$ \\
   \hline\hline
   $x_0^2+x_1^2+x_2^2$  &4            & 3            &{6} \\
   \hline
   $x_0^2+x_1^2$        &1            &{6}           &{6} \\
   \hline
   $x_0^2-x_1^2$        &{7}          & 3            & 3 \\
   \hline
   $x_0^2$              &4            & 9            & 0 \\
   \hline
   \end{tabular}
   }
   \end{equation}
(Here $\#M$ denotes the cardinality of a set $M$.) Therefore a subset
$b$ of $\Wcal$ is a block if, and only if, one of the following
holds:

{\em Case A:} $b=\Qcal_t(q)$ with $t\in F\setminus\{0\}$ and $U\in
\Qcal_{2t}(q)$. By $\#\Qcal_t(q)>3$ and $\#\Qcal_{2t}(q)>0$ the last
two lines of (\ref{planar-tabelle}) can be ruled out. So either
$\Qcal_0(q)$ is a conic ($x_0^2+x_1^2+x_2^2=0$) and $b$ is the set of
its six external points, whereas $U$ is internal, i.e., it does not
lie on a tangent; or $\Qcal_0(q)$ is a singleton ($x_0^2+x_1^2=0$)
and $b$ is the symmetric difference $(r\cup s)\setminus (r\cap s)$ of
two distinct lines with $U\notin r\cup s$ and $\Qcal_0(q)=r\cap s$.

{\em Case B:} $b=\Qcal_0(q)\setminus\{U\}$ with $q\neq 0$ and
$U\in\Qcal_0(q)$. We infer from $\#\Qcal_0(q)>4$ and
(\ref{planar-tabelle}) that $\Qcal_0(q)$ is a pair of lines
($x_0^2-x_1^2=0$), say $\Qcal_0(q)=g\cup h$. Therefore $b=(g\cup
h)\setminus\{U\}$, where $U\in g\cup h$.

Very loosely speaking, a block is either one ``side'' of a quadric
with $U$ being on the ``other side'', or it is the set of all points
in $\Wcal$ of a quadric containing $U$.

   \begin{theo}
The incidence structure $(\Wcal,\Bcal,\in)$ described in Definition
\ref{planar-bloecke} is a $5$--$(12,6,1)$ design.
   \end{theo}

   \proof
By (\ref{point-W}), there are $12$ points in $\Wcal$ and, from our
previous discussion, all blocks have exactly $6$ elements.

In the sequel let $(u_0,u_1,u_2):=(1,0,0)$. So, in terms of
coordinates, an equation of a block takes the form
   \begin{equation} \label{planar-gleichung}
   \sum_{0\leq i\leq j\leq 2} a_{ij}x_i x_j = 2a_{00},
   \mbox{ with at least one } a_{ij}\neq 0.
   \end{equation}
Suppose that we are given a $5$--set $\Dcal=\{F(d_{0k},d_{1k},d_{2k})
\mid k\in\{0,1,\ldots,4\}\}$ contained in $\Wcal$. In order to obtain
all blocks through $\Dcal$ we have to find the non--zero solutions of
the linear homogeneous system
      \begin{equation}\label{planar-LGS}
      \sum_{0\leq i\leq j\leq 2} a_{ij}d_{ik} d_{jk} = 2a_{00},
      \quad k\in\{0,1,\ldots,4\}.
      \end{equation}
This is a system of $5$ equations in $6$ unknowns $a_{ij}$, whence a
non--zero solution exists, i.e., there is at least one block
containing $\Dcal$. In order to show its uniqueness, we have to
distinguish two cases:

{\em Case A:} Each solution with $a_{00}=0$ is trivial. Consequently,
   \begin{equation}
   \det\left(
   \begin{array}{ccccc}
   d_{00}d_{10} & d_{00}d_{20} & d_{10}^2 & d_{10}d_{20} & d_{20}^2 \\
   d_{01}d_{11} & d_{01}d_{21} & d_{11}^2 & d_{11}d_{21} & d_{21}^2 \\
   \multicolumn{4}{c}\dotfill\\
   d_{04}d_{14} & d_{04}d_{24} & d_{14}^2 & d_{14}d_{24} & d_{24}^2
   \end{array}
   \right)
   \neq 0,
   \end{equation}
so that all solutions of (\ref{planar-LGS}) form a one--dimensional
subspace of $F^6$, as required.

{\em Case B:} There is a non--trivial solution $(\qu
a_{00},\ldots,\qu a_{22})\in F^6$ with $\qu a_{00} = 0$. The numbers
$\qu a_{ij}\in F$ determine a non--zero quadratic form \abb{\qu
q}{F^3}{F} and a quadric $\Qcal_0(\qu q)$ containing
$\Dcal\cup\{U\}$. By (\ref{planar-tabelle}), $\Qcal_0(\qu q)$ is a
pair of lines, say $\qu g\cup \qu h$. So
      \begin{equation}
      \qu b := (\qu g\cup \qu h) \setminus \{U\}
      \end{equation}
is one block through $\Dcal$.

Conversely, let $b$ be a block passing through $\Dcal$. There is no
$5$--arc in \PG23, so that at least three points of $\Dcal$ are on a
line, say $\qu g$. There are three possibilities:

1. $b$ stems from a pair of lines $g\cup h$, i.e., $b=(g\cup h)
\setminus \{U\}$. Then the quadrics $g\cup h$ and $\qu g\cup \qu h$
have six common points, whence they are identical.

2. $b$ is the set of external points of a conic $c$. No line contains
four external points of a conic (cf.\ \cite[p.\ 178]{hirs98}). So
$\qu g$ is a tangent of $c$ and $\#(\qu g\cap\Dcal)=3$. We infer
$U\notin\qu g$, since there are no internal points on a tangent. By
$\#\Dcal=5$, $\qu h\setminus\qu g$ contains two distinct external
points and the internal point $U$. Hence $\qu h$ is an exterior line
carrying two distinct internal points. This implies that $\qu g$ and
$\qu h$ meet at an internal point which contradicts $\qu g$ being a
tangent.

3. There are two distinct lines $r$, $s$, with $b=(r\cup s)\setminus
(r\cap s)$ and $U\notin r\cup s$. W.l.o.g.\ let $\#(r\cap\Dcal)=3$,
so that $\#(s\cap\Dcal)=2$. The quadric $\qu g\cup \qu h$ contains
three distinct points of $r$, whence $r\subset \qu g\cup \qu h$.
Similarly, it follows now that $s\subset \qu g\cup \qu h$. Hence
$U\in \qu g\cup \qu h = r\cup s$, an absurdity.

Thus obviously $b=\qu b$.
\qed

   \begin{remark}
   {\em
Up to isomorphism, the Witt design $W_{12}$ is the only
$5$--$(12,6,1)$ design \cite[Chapter IV, \S 2]{beth-jung-lenz85}. The
stabilizer of $U$ in the collineation group of \PG23 yields a
subgroup of the automorphism group of $W_{12}$, i.e., the Mathieu
group $M_{12}$.
   }
   \end{remark}
   \begin{remark}
   {\em
If $\{A,B,C,U\}=:g$ is a line of \PG23, then the three--fold derived
design $(W_{12})_{A,B,C}$ is an affine plane of order $3$. It is
immediate from the definition of blocks that this is just the affine
plane $\Acal$ which arises from \PG23 by removing the line $g$.

Each affinity $\alpha$ of $\Acal$ extends, on one hand, to a unique
collineation $\kappa$ of \PG23 and, on the other hand, to a unique
automorphism $\beta$ of $W_{12}$.

For each $X\in g\setminus \{U\}$ there is a unique elliptic (i.e.,
fixed--point free) involution $\gamma_X$ of $g$ which interchanges
$X$ with $U$. We mention without proof that
   \begin{equation}
   X^\beta = U^{\kappa\inv\gamma_X\kappa}.
   \end{equation}
   Thus $X^\kappa$ and $X^\beta$ need not coincide. Cf.\ also
   \cite[Remark 6]{havl98}.

The discussion of other derivations of $W_{12}$ in terms of the
present planar model is left to the reader.
   }
   \end{remark}

\vspace{1em}
\noindent Hans Havlicek, Abteilung f\"ur Lineare Algebra
und Geometrie, Technische Universit\"at, Wiedner Hauptstra{\ss}e
8--10, A--1040 Wien, Austria.\\ EMAIL: {\tt
havlicek@geometrie.tuwien.ac.at}
\end{document}